\documentclass[10pt,twoside]{article}   
\usepackage{amsmath}
\usepackage{amsfonts}

\usepackage{ecphead}                    

\newtheorem{thm}{Theorem}
\newtheorem{lemma}[thm]{Lemma}

\newtheorem{prop}[thm]{Proposition}
\newtheorem{rem}[thm]{Remark}

\newcommand{\ee}{\end{equation}}
\newcommand{\be}[1]{\begin{equation}\label{#1}}

\newenvironment{proof} {\noindent{\sc Proof. }}{{\hfill
$\Box$}\par\vskip2\parsep}

\usepackage[british]{babel}
\usepackage{amssymb,amscd}
\usepackage[T1]{fontenc}
\usepackage[latin1]{inputenc}
\usepackage{latexsym}
\usepackage{enumerate}
\usepackage{url}
\usepackage{amsbsy}
\usepackage{xspace}

\newcommand{\ie}{{\it i.e. }}
\newcommand{\impt}[1]{{\em #1}}

\newcommand{\nosPCA}{\ensuremath{{\mathcal C}_0}\xspace}

\renewcommand{\emptyset}{\varnothing}

\renewcommand{\leq}{\leqslant}
\renewcommand{\geq}{\geqslant}

\newcommand{\miniop}[3]{%
\renewcommand{\arraystretch}{0.6}
\begin{array}{c}
{\scriptstyle #1}\\
#2\\
{\scriptstyle #3}
\end{array}
\renewcommand{\arraystretch}{1}}

\newcommand{\card}[1]{\ensuremath{\# #1}}

\newcommand{\eqdef}{\triangleq}

\newcommand{\argth}{\ensuremath{\textrm{Argth}}}

\renewcommand{\L}{\Lambda}
\newcommand{\s}{\ensuremath{\sigma}}
\newcommand{\e}{\ensuremath{\eta}}
\renewcommand{\r}{\ensuremath{\rho}}
\renewcommand{\o}{\ensuremath{\omega}}
\renewcommand{\t}{\ensuremath{\tau}}
\renewcommand{\d}{\ensuremath{\delta}}
\renewcommand{\b}{\ensuremath{\beta}}
\newcommand{\z}{\ensuremath{\zeta}}

\renewcommand{\epsilon}{\varepsilon}

\newcommand{\Zd}{\mathbb Z^d}
\newcommand{\Z}{\mathbb Z}
\newcommand{\N}{\mathbb N}

\newcommand{\tq}{ : \ }

\newcommand{\ind}[2]{1 \hspace{-0.9ex} 1_{#1} {#2}}

\newcommand{\spin}{S{}}

\newcommand{\cvinf}{\ensuremath{\spin^{\Zd}}}

\newcommand{\cvinfLx}[1]{\ensuremath{\spin^{#1}}}

\newcommand{\SpinConfigSpace}{\ensuremath{\{-1,+1\}^{\Zd}}}

\newcommand{\abs}[1]{\ensuremath{ \Big|  #1  \Big| }}

\newcommand{\norme}[1]{\ensuremath{ \| #1 \|_{_{_1}} }}

\newcommand{\tvert}{\vert \hspace{-0.8ex} \parallel}

\newcommand{\ntrois}[1]{\tvert #1 \ \tvert}

\newcommand{\boule}[2]{\ensuremath{ {\mathcal B} (#1,#2)}}

\newcommand{\BL}{\ensuremath{ {{\mathcal B}({L})}}}

\newcommand{\BLR}{\ensuremath{ {{\mathcal B}({L-R})}}}

\newcommand{\BnR}{\ensuremath{ {{\mathcal B}({nR})}}}

\newcommand{\nusup}{\ensuremath{{\nu^\maximalConf}}}
\newcommand{\nuinf}{\ensuremath{{\nu^\minimalConf}}}

\renewcommand{\preceq}{\preccurlyeq}
\renewcommand{\succeq}{\succcurlyeq}

\newcommand{\coup}[1]{\mathbf{I}\hspace{-0.5ex}\mathbf{P}^{#1}}
\newcommand{\coupling}[2]{\coup{} \Big( #1 \Big| #2 \Big)}

\newcommand{\opreceq}{\circledast}

\newcommand{\coupCond}[1]{\mathbf{l}\hspace{-1.1ex}\mathbf{Q}^{\maximalConf,\maximalConf} {( #1 ) }}

\newcommand{\W}{\ensuremath{{\mathcal W}}}
\newcommand{\Zp}{\ensuremath{{\mathcal Z}}}

\newcommand{\K}{\ensuremath{{\mathcal K}}}

\newcommand{\G}[1]{\ensuremath{{\mathcal G}(#1)} }

\renewcommand{\S}{\ensuremath{{\mathcal S}} }
\newcommand{\Ss}{\ensuremath{{\mathcal S}_{s}} }

\newcommand{\R}{\ensuremath{{\mathcal R}} }
\newcommand{\Rs}{\ensuremath{{\mathcal R}_{s}} }

\newcommand{\di}[1]{\ensuremath{\partial_i #1}}

\newcommand{\minimal}{\ensuremath{{-}}}
\newcommand{\maximal}{\ensuremath{{+}}}
\newcommand{\minimalConf}{\ensuremath{{\pmb{-}}}}
\newcommand{\maximalConf}{\ensuremath{{\pmb{+}}}}

\newcommand{\assump}{\ensuremath{(\mathcal A)}\xspace}
\newcommand{\WM}{\ensuremath{(\mathcal W \mathcal M)}\xspace}

\newcommand{\+}{{\bf +1}}

\newcommand{\Var}[2]{\ensuremath{\Delta_{#2}(#1)}}

\newcommand{\constante}{\ensuremath{\kappa}\xspace}

\begin{document}


\input{paper13.head.sty} 


{\bf REPRINT,  DOI: 10.1214/ECP.v9-1116}

\section{Introduction}

  The main feature of
  Probabilistic Cellular Automata dynamics  (usually abbreviated in PCA)
   is the parallel, or
  synchronous, evolution of
  all interacting elementary components.
    They are precisely discrete-time Markov
  chains on a product space $S^{\L}$ (\impt{configuration space})
  whose transition probability is a \impt{product measure}. In this paper,
$\spin$ (\impt{spin space}) is assumed to be a finite %
set with total order
denoted by $\leq$ and $\L$
(set of \impt{sites}) a subset, finite or infinite, of $\Zd$.
The fact that the transition probability kernel $P(d\s|\s')$ %
($\s,\s' \in
S^{\L}$) is a product measure means that all spins $\{\s_{k}: k
\in \L\}$ are simultaneously and independently updated.
This transition mechanism differs from the
one in the most common Gibbs samplers, where only
one site is updated at each time step.
In opposition to these dynamics with sequential updating, it is simple to
define PCA's on the infinite set $\cvinf$ without passing to
continuous time.
\medskip

The main purpose of this article is to study relation between different types of conditions
which insure the
 fastest convergence towards an equilibrium state ($\nu P=\nu$)
of PCA dynamics on $\cvinf$.
 Let us emphasise that the  non-degeneracy hypothesis we will assume implies that the
 asymptotic behaviour
of PCA dynamics on $\cvinfLx{\L}$ where $\L$ is a finite subset of~$\Zd$ (called
\impt{finite %
volume PCA dynamics}) is well-known.
It is a classical result from the theory of finite state space aperiodic
irreducible Markov Chains. Such discrete time processes
admit a unique stationary probability
measure, and are ergodic.
However, if the PCA dynamics is considered on $\cvinf$ (\impt{infinite
  volume dynamics}), some
non-ergodic behaviour may arise (see for instance example~2 %
section~$I\!I\!I$ %
{\it in}~\cite{KV}).
The most famous condition which insures
ergodicity of the PCA dynamics on
$\cvinf$ is due to Dobrushin and Vasershtein's work
(see~\cite{Vasershtein}), and applies in the
high-temperature
regime. Others conditions of ergodicity for general PCA
can be
found in the following
works:~\cite{Ferrari91,MalyshevIgnatiuk88,LMS,%
MaesShlosman91,MalyshevMinlos}.
See for instance Sections~6.1.2 and~6.1.3 in~\cite{ThesePyl} for details.
They all are effective only when some high-temperature condition holds or
in perturbative cases.
\medskip

We will here adopt another approach, partially inspired by
Martinelli and Olivieri's work for a class of continuous time
Interacting Particle Systems
called Glauber dynamics
(see~\cite{MartinelliOlivieriI}),
and based on a famous statement of Holley about
rate of convergence~(\cite{Holley}).
We introduce a condition~\assump which means
the \impt{exponential decay of the influence from the boundary for the %
invariant measure of the system restricted to any finite box},
which will be here proved to be equivalent to the exponentially fast ergodicity~(Theorem~\ref{ergodicite}).
The condition~\assump we use is not a constructive criterion
like the Dobrushin-Vasershtein condition, or
its generalised version
developed in~\cite{MaesShlosman91} and numerically studied in~\cite{deJongMaes}.
But, theoretically, comparison of spatial and time mixing are always %
interesting (cf.~\cite{MartinelliOlivieriI,Dyer_et_al}).
Furthermore we present different examples in which~\assump
is satisfied on a larger domain than Dobrushin-Vasershtein condition, and is moreover optimal for these models.
\medskip

In section~\ref{Main} we state our main results.
The first and more general one (Theorem~\ref{ergodicite}) is the following: convergence towards equilibrium
in the uniform norm with an \impt{exponential rate} is \impt{equivalent} to the condition~\assump .
In other words \impt{exponential mixing in space is equivalent to exponential mixing in time}.
It will then be applied to  %
a class of reversible PCA dynamics on~\SpinConfigSpace,
associated in a natural way to a Gibbsian potential $\varphi$. We prove
that the usual weak mixing condition for~$\varphi$ implies the validity
of~\assump, thus the {\it exponential ergodicity} of the
dynamics
towards the unique Gibbs measure associated to $\varphi$ holds
(Theorem~\ref{ErgodClasseReversible}).
For some particular PCA of this class,
we also prove that~\assump
is weaker than the Dobrushin-Vasershtein ergodicity
condition and note that the exponential
ergodicity holds as soon as there is no phase transition.
Our result are then the first optimal ones in this context.
Sections~\ref{SectionProof1} and~\ref{SectionProof2} are respectively %
devoted to the proof of the Theorems and useful Lemmas.

\section{\label{Main}Main results}

Let $P$ denotes a PCA dynamics on
$\cvinf$. This means a Markov Chain on \cvinf whose
transition probability kernel $P$ verifies
 for all \impt{configuration} $\eta \in
\cvinf$, \mbox{$\s=(\s_k)_{k \in \Zd} \in \cvinf$},
$ P(\ d\s \ | \ \e \ )=\miniop{}{\otimes}{k \in \Zd} p_k(\ d\s_k \ | %
\ \e \ ) $,
where for all \impt{site} $k \in \Zd$, for all $\eta$,
$p_k(\ .\ |\e)$ is a probability measure on~$\spin$, called
\impt{updating rule}.
For any subset $\Delta$ of $\Zd$, and for all configurations
$\s$ and $\e$ of $\cvinf$, the configuration
$\s_\Delta \e_{\Delta^c}$ is defined by $\s_k$ if $k \in \Delta$,
else $\e_k$.
Let the notation $\s_\Delta$ design $(\s_k)_{k \in \Delta}$ too.
Let $\L$ be a finite subset of $\Zd$ (denoted by \mbox{$\L \Subset
\Zd$}). We call \impt{finite volume PCA dynamics} with \impt{boundary
condition} $\tau$ ($\tau \in \cvinf$ or $\tau \in
\cvinfLx{\L^c}$),
the Markov Chain on
$\cvinfLx{\L}$ whose transition probability $P_\L^\t$
is defined by:
$ P_\L^\t( d\s_\L \ | \ \e_\L\ )=\miniop{}{\otimes}{k \in \L} %
p_k(\ d\s_k \ | \ \e_\L\t_{\L^c} \ ) .$
It may be identified with the following infinite volume %
PCA dynamics on~$\cvinf$:
$ P_\L^\t( d\s \ | \ \e_\L\ )=\miniop{}{\otimes}{k \in \L} %
p_k(\ d\s_k \ | \ \e_\L\t_{\L^c} \ ) \otimes  \d_{\t_{\L^c }}%
(d\s_{\L^c})$.
Let $\nu_\L^\t$ denote the stationary
measure associated to the finite volume
dynamics $P_\L^\t$.
For $\nu$ probability measure on $\cvinf$
(equipped with the Borel $\s$-field associated to the product
topology),
$\nu P$ refers to
$\nu P (d\s)=\int P(d\s |  \e ) \nu(d\e)$.
Recursively $\nu P^{(n)}=(\nu %
P^{(n-1)}) P$.
For each function $f$ on $\cvinf$, $P(f)$ is the function defined by
$P(f)(\e)=\int f(\s) P(d\s|\e)$.
All the measures considered
in this paper are probability measures.

\medskip

PCA dynamics considered here are assumed to be \impt{non degenerate}:
$\forall k \in \Zd,\ \forall \e \in \cvinf$, \mbox{$\forall s \in \spin$},
$ p_k(\ s \ | \ \e \ ) >0 $; they are also \impt{local}, which means:
$\forall k \in \Zd, \exists \ V_k \Subset \Zd, %
p_k(\ .\ |\e)=p_k(\ .\ |\e_{V_k}) $ and
they are also translation invariant:
$ \forall k \in \Zd,\ \forall s \in \spin,\ \forall \e \in \cvinf,\quad %
p_k(\ s \ | \ \e \ )= %
p_0(\ s \ | \ \theta_{-k}\e \ )$,
where  $\theta_{k_0}(\s)$ defines the
translation of a configuration %
$\s$ of $\cvinf$ with
\mbox{$\theta_{k_0}(\s) =$} \mbox{$ (\s_{k-k_0})_{k \in \Zd}$}.
\medskip

 Attractivity of PCA dynamics is moreover assumed here:
 One can order two configurations by defining $\s \preceq \e$ if $\forall k \in \L, \s_k \leq \e_k$.
 A real function $f$ on $\cvinfLx{\L}$ will then be said to be increasing if $\s \preceq \e$ %
implies $f(\s) \leq f(\e)$. Thus two probability measures $\nu_1$ and $\nu_2$ satisfy %
 the stochastic ordering $\nu_1 \preceq \nu_2$ if,
 for all increasing functions~$f$ on $\cvinfLx{\L}$, $\nu_1(f) \leq \nu_2(f)$, with %
the notation $\nu_i(f)=\int f(\s) \nu_i (d\s)$.
 As Markov chain, a PCA dynamics $P$ on $\cvinfLx{\L}$ ($\L \subset \Zd$) is \impt{attractive} if
 for all increasing function $f$, $P(f)$ is still increasing.
 Let us define too, for $s \in \spin,\s \in \cvinfLx{\L}$,
 the function $G_k(s,\s)$ by:
 \be{ContFonctRepart}
 G_k(s,\ . \ )= \sum_{s'\geq s} p_k(s'|\ .\  ).
 \ee
 Recall that a PCA dynamics is attractive if, and only if,
 for all $k$ in $\L$, and all value  $s \in \spin$, the function
 $G_k(s,.)$ is increasing (in $\s$).
 \medskip

 A real valued function~$f$  on \cvinf is said \impt{local}
if
$ \exists \L_f \Subset \Zd,\ \forall \s \in \cvinf, \ %
f(\s)=f(\s_{\L_f}) $.
We define, for each $f$ continuous function on the compact \cvinf and
for all $k$ in $\Zd$,  %
$$  \Var{k}{f}=\sup \Big\{ \abs{f(\s)-f(\e)} \ : \ (\s,\e)\in (\cvinf)^2, %
\s_{\{ k  \}^c} \equiv \e_{\{ k  \}^c} \Big\} ,$$
and the semi-norm %
\mbox{$ \ntrois{f} = \sum_{k \in \Zd} {\Var{k}{f}}$}.
For $L$ integer, $\BL$ is the ball
 $\boule{0}{L}$
with respect to the norm \mbox{$\norme{k} = \sum_{i=1}^{d} | k_i |$},
\mbox{$k=(k_1,k_2,\hdots,k_d) \in \Zd$}.

\begin{thm}
\label{ergodicite}
Let $\spin$ be a totally ordered finite set with maximal
(resp.~minimal) element denoted by~\maximal (resp.~\minimal).
$\maximalConf$ (resp. $\minimalConf$) denotes configurations equal to $\maximal$ (resp. $\minimal$) in all sites.
Let $P$ be an attractive, translation invariant, non degenerate, local
PCA
dynamics on \cvinf.
Let $\nu_\BL^{\maximalConf}$ (resp.~$\nu_\BL^{\minimalConf}$) be the stationary measure
of $P_\BL^\maximalConf$ (resp.~$P_\BL^\minimalConf$).
The following spatial mixing condition:
\mbox{$\exists C>0, \ \exists M > 0, \ \exists L_1 \in \mathbb N^*, \forall L \in \mathbb N^*, L \geq L_1,$}
$$ \label{ExpControl} %
\int \ \s_0 \ d \nu_{\BL}^{\maximalConf} - %
\int \ \s_0 \ d \nu_{\BL}^{\minimalConf}  \leq C e^{-ML}
\qquad \qquad \qquad \assump  $$
is equivalent to the convergence
of the dynamics~$P$ towards the unique
equilibrium state~$\nu$ with exponential rate:
$ \exists \lambda >0$,  $\exists n_1$, $\forall n \geq n_1$, %
$\forall f \textrm{ local function on } \cvinf :$
\be{ErgodExp}
\sup_{\s} \abs{ %
\d_\s P^{(n)} (f) - \nu(f) } \leq %
2 \ntrois{f} e^{-\lambda n} . \ee
\end{thm}

In order to better interpret the meaning of condition~\assump
and the relevance of Theorem~\ref{ergodicite}, we then apply
it to a wide class %
of reversible PCA dynamics on $\SpinConfigSpace$.
First, let us recall some known facts about \impt{reversible PCA dynamics}
(that is to say PCA dynamics
whose set of reversible measures~$\R$ is not empty).
The study of the qualitative nature  of their equilibrium states
as Gibbs measures
was initiated by
Kozlov and Vasilyev (see~\cite{KV,Vasilyev}).
Gibbs measures
with respect to some
dynamics' naturally associated potential, are indeed natural candidates as
stationary states.
In~\cite{DPLR,ThesePyl},
precise relations were established
between the sets of stationary measures, reversible measures and
some Gibbs measures
(see Proposition 3.3 in~\cite{DPLR}).
Moreover, unlike what is done (or expected to hold) for continuous time
Interacting Particle Systems like Glauber dynamics or gradient
diffusions, it is shown
that Gibbs measures may be
non stationary for PCA's dynamics, which is a characteristic
manifestation of the discrete time case.
\medskip

Assume until the end of this section and in section~\ref{SectionProof2} that \mbox{$\spin=\{-1,+1\}$}.
We call
class~\nosPCA
the family of PCA dynamics on~$\SpinConfigSpace$ %
whose updating rule
 $(p_k)_{k \in \Zd}$ is given by: $\forall k \in \Zd$, $\forall \e \in \cvinf$,  %
$ \forall s \in \spin$
\be{FormeLocaleNosPCAPrimaire}
 p_k(s \ | \ \e)=\frac{1}{2} \Big(1+ s \tanh ( %
\b \sum_{k'\in \Zd} \K(k'-k) \e_{k'}  ) \Big),
\ee
where~$\b$ is a positive real parameter and \mbox{$\K\ : \ \mathbb Z^d \rightarrow \mathbb R$}
is an interaction function between sites
which is symmetric
and has finite range $R>0$
(\ie
for all $k$ of $\Zd$ such that
$ \norme{k} > R$ then $\K(k)=0$).
Remark  that $\b = 0$ is the independent case (sites don't interact), %
and that when $\b$ increases, the dynamics becomes less and less random.
So $\b$ may be thought as a kind of inverse temperature parameter.
See subsection~4.1.1 in~\cite{ThesePyl} for the generality of the  class
\nosPCA among reversible PCA dynamics on~$\SpinConfigSpace$.
Due to their definition, PCA dynamics in \nosPCA are local, %
translation invariant, non degenerate.
It is known (see~\cite{KV,DPLR}) that
any PCA dynamics~$P$ in~\nosPCA
admits at least one reversible measure which is a Gibbs measure associated to %
the following translation invariant multibody potential~$\varphi$:
 \be{DefPotentiel}
 \begin{array}{lll}
     \varphi_{U_{k}}(\s_{U_k}) & = & - \log \cosh \left( \b  \sum_{j} \K(k-j)
     \s_{j} \right) \qquad \textrm{where } U_k =\{ j : \K(k-j) \neq 0 \}\\
     \varphi_{\L}(\s_{\L}) & = & 0 \ \ \text{otherwise}.
 \end{array}
 \ee
Moreover Proposition~3.3 in~\cite{DPLR} stated the precise relations
$\R=\S \cap \G{\varphi}$ and
$\Rs = \Ss$,
where~$\S$ (resp.~$\R$) denotes the set of $P$-stationary
(resp.~$P$-reversible) measures, $\Ss$ and $\Rs$ their respective
space-translation invariant measures' parts, and $\G{\varphi}$ the set of
Gibbs measures on $\cvinf$ %
associated to the potential~$\varphi$.

One also checks that such a PCA dynamics $P$ is attractive, %
if and only if function $\K(.)$ is non-negative %
(see Property~4.1.2 in \cite{ThesePyl}).
From now on, let us assume that $\K$ is non negative.
\medskip

Mixing conditions for a potential~$\varphi$ define
different regions in the domain of absence of phase transition for the associated Gibbs measures.
Strong mixing conditions are
usually related to the domain where Dobrushin's uniqueness holds,
and weak mixing conditions are expected to be valid in the main part
of the uniqueness domain: See~\cite{MartinelliOlivieriI}
for a review on these conditions.
Here, we call weak mixing condition
 for the potential $\varphi$, the condition:
\mbox{$\exists C>0,\ \exists M>0,\ \forall L \geq 2,$}
$$ \int \ \s_0 \ \mu(d\s_\BL | \s_{\BL^c}= +1 ) - %
\int \ \s_0 \ \mu(d\s_\BL | \s_{\BL^c}= -1 ) \leq %
C e^{-ML} \qquad \WM $$
where $\mu$ is the unique Gibbs measure associated to~$\varphi$.
For ferromagnetic potentials, %
it is indeed the equivalent form of more general weak mixing condition.

\begin{thm} \label{ErgodClasseReversible}
Let $P$ be an attractive %
PCA dynamics on \SpinConfigSpace of the class \nosPCA %
defined by~(\ref{FormeLocaleNosPCAPrimaire}),  let $\varphi$
denote the potential canonically associated
defined in~(\ref{DefPotentiel}), and $\G{\varphi}$ the set of Gibbs measures
w.r.t $\varphi$.
\begin{itemize}
\item If there is phase transition (\ie \mbox{$\card{\G{\varphi}} >1$})
then the dynamics~$P$ is
  non-ergodic.
\item Otherwise, when there is no phase transition
(\ie  \mbox{$\G{\varphi}=\{ \mu \}$})
the dynamics~$P$
  is ergodic towards the unique Gibbs measure
  $\mu$. \\
Moreover if we
  assume the potential~$\varphi$ satisfies the weak mixing condition~\WM, then the convergence towards
  $\mu$ holds with exponential rate.
\end{itemize}
\end{thm}

In~\cite{DPLR}, we established that, for nearest neighbour interaction
function $\K$, phase transition holds for $\b$ large.
For instance, when
$d=2$, let $P_J$~be the PCA dynamics of the class~\nosPCA obtained
taking:
$ \K(\pm e_1)=\K(\pm e_2)=J>0,\ \K(k)=0 $ otherwise,
where $(e_1,e_2)$ is a basis of $\mathbb R^2$ and $J$ a positive constant.
The canonically associated potential~$\varphi_J$ ({\it cf.}~(\ref{DefPotentiel}) )
is the following four-body potential:
$ \varphi_{J,V_{k}}(\s_{V_k})  = - %
\log \cosh ( \b J \sum_{j \in U_k} \s_j )$ %
where \mbox{$U_k=\{ k-e_1,k+e_1,k-e_2,k+e_2 \}$}.
From Theorem~\ref{ErgodClasseReversible}
we conclude here that for $\b$ large, the PCA $P_J$ is non-ergodic since it
has at least two different stationary states $\nu^-$ and $\nu^+$.

\medskip

Let us now discuss how large is the domain where condition~\WM holds.
One conjectures Weak Mixing condition for Gibbs measure is
valid up to  the critical temperature, that is, as soon as there
is no phase transition. In that sense, our main result would give
ergodicity with exponential rate on a much larger region as
the region where the Dobrushin-Vasershtein criterion holds.
In fact, let us mention
the reference~\cite{Higuchi},  where, using percolation techniques,
it is proved that in
dimension $d=2$, for a ferromagnetic nearest neighbour Ising model
without extremal magnetic field, the associated Gibbs measure is weak mixing
as soon as it is unique (\ie $\forall \b, \b < \b_c$).
In order to precise this assertion, let us consider
the dynamics~$P_J$.
A projection argument relates the
potential~$\varphi_J$ associated to~$P_J$
with the usual Ising ferromagnetic pair potential with intensity coefficient~$J$
(see~\cite{Vasilyev}). Due to Higuchi's result, we know that
the Gibbs state associated to this potential~$\varphi_J$
is weak mixing as soon as there is no phase transition, which happens for
$\b$ lower than the critical value $\b_c$, which coincides with
the Ising critical inverse temperature
\mbox{$\b_c=\frac{\log (1+\sqrt{2})}{2J}$}.
In other words, we obtain that the PCA dynamics $P_J$ is %
ergodic with exponential rate for $\b < \b_c$ and non-ergodic for $\b>\b_c$.
Taking $J=1$, \mbox{$\b_c \simeq 0.441$}; since Dobrushin-Vasershtein criterion
applies only for
\mbox{$\b<\frac{1}{2} \argth (\frac{1}{2}) \simeq 0.275 $}
({\it cf.}~part~6.1.2 {\it in}~\cite{ThesePyl}), ours is better.

\section{\label{sectionErgodicite}\label{SectionProof1}Proof of the Theorem~\ref{ergodicite} }

The proof of Theorem~\ref{ergodicite} is based on the existence of some coupling of %
PCA dynamics preserving the stochastic ordering.
Let  $(P^1,P^2,\hdots,P^N)$ be an increasing  $N$-uple of PCA dynamics
which means PCA related by the following monotonicity property
$  \forall k \in \Zd$,
$\forall \z^1 \preceq \z^2 \preceq \hdots \preceq \z^N \in \cvinf, \forall s \in \spin $,
$ G_k^1(s \ | \ \z^1 ) \leq G_k^2(s \ | \ \z^2 ) \leq \hdots %
\leq G_k^N(s \ | \ \z^N ) $ where $G^i$ is the function associated to~$P^i$ by~(\ref{ContFonctRepart}).
There exists~(cf.~\cite{LouisCoupling})   a monotone synchronous coupling on $(\cvinf)^N$
denoted by %
\mbox{$P^1 \opreceq P^2 \opreceq \hdots \opreceq P^N$}
with the following property:
for all initial configuration
\mbox{$ \s^1 \preceq \s^2 \preceq \hdots \preceq \s^N $} %
 and for all times~$n$,
$$ P^1 \opreceq \hdots \opreceq P^N \ \big( \ %
\o^1(n) \preceq \hdots \preceq \o^N(n) \ \big| \ (\o^1,\hdots,\o^N)(0)= %
(\s^1,\hdots,\s^N) \ \big) =1 . $$
Such a coupling %
will be called \impt{increasing synchronous coupling}.
The notation $\coup{}$
denotes the coupling $P \opreceq P \opreceq \hdots \opreceq P$
of $N$~times the same PCA dynamics $P$, where $N$ will be a finite
large enough number.
\medskip

This coupling allows us to develop some monotonicity argument and to
state the following result, whose proof is in~\cite{LouisCoupling}:
\begin{prop} \label{ResultatsDeMonotonie}
The measure
 $\nu_\L^{\maximalConf}$ (resp. $\nu_\L^{\minimalConf}$)
is  the maximal (resp. minimal) measure of the set
\mbox{$\{ \nu_\L^\t \tq \t \in \cvinfLx{\L^c} \}$} of stationary
measures associated to the PCA dynamics~$P_\L^\t$ on the fixed finite volume~$\L$ %
and with boundary condition~$\t$.
Let $\nusup$ and $\nuinf$ denote the maximal and the minimal elements %
of the
set~$\S$ of stationary measures associated to the PCA dynamics~$P$.

Following relations hold:
\be{LimTempLimSpatPlus}
\nusup= \lim_{L \to  \infty} \nu_\BL^{\maximalConf} \otimes %
\delta_{(\maximalConf)_{\BL^c}}=\lim_{n \to \infty} \delta_{\maximalConf} P^{(n)}
\ee %
\be{LimTempLimSpatMinus} %
\nuinf= \lim_{L \to  \infty} \nu_\BL^{\minimalConf} \otimes %
\delta_{(\minimalConf)_{\BL^c}}=\lim_{n \to \infty} \delta_{\minimalConf} P^{(n)}.
\ee
In particular, $P$ admits a unique stationary measure~$\nu$ if and
only if %
\mbox{$ \nuinf=\nusup $}.
\end{prop}

Note that $P^{(n)}$ denotes $P \circ P \circ \hdots \circ P$, and so %
is for instance $\delta_{\maximalConf} P^{(n)}$ the law
at time~$n$ of the Markov Chain with transition kernel~$P$ and initial %
distribution~$\delta_{\maximalConf}$.

\begin{rem} \label{discreteTime}
{\rm Note the following range of dependence w.r.t. the past for local PCA.
Let us define
$\overline{\L} = \cup_{k \in \L} V_k =\overline{\L}^{(1)},$
and
$\overline{\L}^{(n)}%
  = \cup_{k \in {\overline{\L}}^{(n-1)}} V_k $.
Then: $\forall n, \forall \L \Subset \Zd$, $ \forall (\s,\e) \in (\cvinf)^2$ with %
$ \s_{\overline{\L}^{(n)}} \equiv \e_{\overline{\L}^{(n)}}$,
$ 
\coupling{\o^1_\L(n) %
  \equiv \o^2_\L(n)}{(\o^1,\o^2)(0)=(\s,\e)}=1 $.
}
\end{rem}

\begin{proof} {((\ref{ErgodExp}) implies~\assump in Theorem~\ref{ergodicite})} \\
It uses a usual strategy and takes advantage of the coupling~$P \opreceq P_\BL^{\maximalConf}$.
Let $L$ be a fixed integer, larger than $L_1=n_1$ where $n_1$ is defined in~(\ref{ErgodExp}).
Using the relation (stated in \cite{LouisCoupling})
\be{IntermedEncadrementNu}
\nu_\BL^{\minimalConf} \otimes \delta_{(\minimalConf)_{\BL^c}} \preceq \nu %
\preceq \nu_\BL^{\maximalConf} \otimes \delta_{(\maximalConf)_{\BL^c}};
\ee
the positivity of each following term is stated. We have:
$$ 0 \leq  
\int \ \s_0 \ d \nu_{\BL}^{\maximalConf} - %
\int \ \s_0 \ d \nu_{\BL}^{\minimalConf} = %
 \Big( 
\int \ \s_0 \ d \nu_{\BL}^{\maximalConf} - %
\int \ \s_0 \ d \nu \Big) + %
 \Big( 
\int \ \s_0 \ d \nu - %
\int \ \s_0 \ d \nu_{\BL}^{\minimalConf} \Big) ,$$
and we will state that each part is lower than $2 \ntrois{f_0} e^{- \lambda L}$ (where $f_0(\s)=\s_0$).
We only give the proof for \mbox{$ 
\int \ \s_0 \ d \nu_{\BL}^{\maximalConf}  - \int \ \s_0 \ d \nu $}
since the proof
for the minimal~$\minimalConf$ boundary condition is analogous.
For any $n \in \mathbb N^*$,
\begin{eqnarray*} {  \nu_{\BL}^{\maximalConf} (\s_0) - %
\nu (\s_0)  = \Big( \nu_{\BL}^{\maximalConf} (\s_0)  %
- \delta_{\maximalConf} {P_{\BL}^{\maximalConf}}^{(n)} (f_0)  \Big) %
 + \Big(   \delta_\maximalConf {P_\BL^\maximalConf}^{(n)} (f_0) - %
\delta_\maximalConf P^{(n)} (f_0) \Big) +} \\
  { \hspace{8cm} \Big(   \delta_\maximalConf P^{(n)} (f_0) - \nu (\s_0) \Big) .}
\end{eqnarray*}
 Using the monotonicity of \mbox{$ P_\BL^{\maximalConf} \opreceq P_\BL^{\maximalConf} $}
the first term is non positive.
 Using the assumption~(\ref{ErgodExp}) the third term is bounded from %
above by \mbox{$ 2 \ntrois{f_0} e^{- \lambda n}$} \mbox{($\forall n \geq n_1$)}.
 Choose now \mbox{$n=L$}. Rewrite the second term as %
$ \coupCond{ \o_0^2 (n) - \o_0^1(n)}$ where
$$ \coupCond{\ . \ } = P \opreceq P_\BL^{\maximalConf} \big( \ . \ | %
(\o^1,\o^2)(0)=(\maximalConf,\maximalConf) \big) . $$
Using Lemma~\ref{NeqEsp}, we bound the second term from above %
with \mbox{$\kappa ' \  \coupCond{ \o_0^2 (n) \neq  \o_0^1(n)} $}.
According to the construction of the coupling and using Remark~\ref{discreteTime}, note
that with respect to $\coupCond{.}$, \mbox{$\o_0^2 (n) \neq  \o_0^1(n)$} is possible
only if it exists a previous time~$n'$ \mbox{($0 < n' < n$)} and %
a site~$k$ in \mbox{$ \BL^c \cap \overline{\{ 0 \} }^{(n')} $} such that
\mbox{$ \o_k^2 (n') =\o_k^1 (n') \neq \maximalConf$}.
By taking~$n=L$, we have \mbox{$  \overline{\{ 0 \} }^{(n')} \subset \BL $}; so is this
event empty, which ensures \mbox{$\coupCond{\o_0^2(n) \neq \o_0^1(n)} = 0 $}.
Thus is~\assump proved.
\end{proof}

\medskip

\begin{proof} {(\assump implies (\ref{ErgodExp}) in Theorem~\ref{ergodicite})} \\
The most delicate part
is to establish the exponential rate of convergence
towards equilibrium.
Our proof is inspired by Martinelli and Olivieri
proof of
exponential ergodicity for continuous time Glauber dynamics on
$\SpinConfigSpace$ (see~\cite{MartinelliOlivieriI}).
For any time~$n \in \N$, let us define
a coefficient which controls the ergodicity:
\be{DefRho}
\r(n)=\coupling{\o_0^1(n) \neq \o_0^2(n)}{(\o^1,\o^2)(0)%
=(\minimalConf,\maximalConf)} .
\ee
If we assume the exponential bound~\assump,
thanks to forthcoming Lemma~\ref{HWMErgodicite}, we deduce that
\mbox{$\lim_{n \to \infty} \r(n)=0$}.
Reporting assumption \assump in
 the inequality~(\ref{InegRecGeneNouvelle}), we can use %
forthcoming Lemma~\ref{DecrescenzaPolinomiale} to deduce that %
$(\r(n))_{n \in \mathbb N^*}$ converge to $0$ faster than $\frac{1}{n^d}$.
Finally, using inequality~(\ref{InegRho}) and %
Lemma~\ref{VitesseExplicite}, we conclude
that $\r(n)$ converges to $0$ exponentially fast; thus, thanks to %
Lemma~\ref{RhoControlErgod}, conclusion holds.
\end{proof}

{\bf Technical lemmas:} First remark the easy fact:
\begin{lemma}
 \label{NeqEsp}
Let  $(\Omega,{\mathfrak A},{\mathcal P})$ be a probability space, and $Z$
a random variable with values in a finite set $\{z_1 <
\hdots < z_m \}$ of $\mathbb R$, such that
\mbox{${\mathcal P}(Z \geq 0)=1$}.
Then, if
\mbox{$ \constante = \max \{ \frac{1}{z_i}, z_i > 0 , %
1 \leq i \leq m \} $} and $\kappa '=max \{ z_i, 1 \leq i \leq m \} $
(which do not
depend on the law of~$Z$ under~$\mathcal P$) we have:
\mbox{$ {\mathcal P}(Z \neq 0) \leq \constante \ \int Z d \mathcal P $} and
\mbox{$ \int Z d \mathcal P \leq \kappa ' {\mathcal P}(Z \neq 0) $}.
\end{lemma}

Using the monotonicity property of the coupling, the two following Lemmas are easily proved.
\begin{lemma} \label{pptyCoup}
{\tiny $\blacksquare$} $\forall \s,\e \in \cvinf, \ \s \preceq \e$,
$ \coupling{\o_0^1(n) \neq \o_0^2(n)}{(\o^1,\o^2)(0)=(\s,\e)}%
  \leq \r(n) $.

 {\tiny $\blacksquare$} $\forall \L \Subset \Zd$, $\forall n \in \mathbb N$, $\forall \xi \in \cvinf$, \\
$ P_\L^{\minimalConf} \big( \o (n) \in .  \big| %
\o (0)=\xi_\L ({\minimalConf})_{\L^c} \big)  %
  \preceq %
P \big( \o(n)\in .  \big  |   \o(0)=\xi \big) \preceq  %
 P_\L^{\maximalConf} \big( \o(n) \in .  \big|  \o(0)=\xi_\L %
({\maximalConf})_{\L^c} \big)$ . %

{\tiny $\blacksquare$} $ \r(n) \leq P_\L^\minimalConf \opreceq P_\L^\maximalConf %
 ( \o^1_0(n)  \neq \o^2_0(n) \ | %
(\o^1,\o^2)(0)=(\minimalConf,\maximalConf) ) $.
\end{lemma}

\begin{lemma} \label{RhoControlErgod}
The sequence
$(\r(n))_{n \in \mathbb N^*}$ is decreasing, and
$\forall f$, $\forall \s, \e$,
$$ \abs{ %
P(f(\o(n)) | \o(0)=\s ) - P(f(\o(n)) | \o(0)=\e ) } %
 \leq \  2 \ \ntrois{f} \r(n) .$$
Thus, if
\mbox{$\lim_{n \to \infty} \r(n)=0$},
the dynamics $P$ is ergodic, and
\mbox{$ \sup_{\s} \abs{ %
P(f(\o(n)) | \o(0)=\s ) - \nu(f) }$} %
 $ \leq %
 2 \ \ntrois{f} \ \r(n) $,
where $\nu$ denotes the unique stationary measure.
\end{lemma}

Note that due to the monotonicity of $\r(.)$, we can restrict
ourselves to the case $\r(.)>0$.
\begin{lemma}  \label{HWMErgodicite}
$\exists \constante$, $\forall \L \Subset \Zd$, %
$\lim_{n \to \infty} \r (n) \leq \constante %
\big( \int \s_0 \ d\nu_\L^{\maximalConf} - %
\int \s_0 \ d\nu_\L^{\minimalConf}   \big) $.
\end{lemma}

\begin{proof}
Note
$ P_\L^{\minimalConf} \opreceq P_\L^{\maximalConf} %
\Big( \o_0^1(n) \leq \o_0^2(n)) \ \Big| \ %
{(\o^1(0),\o^2(0)) =(\minimalConf,\maximalConf)}
\Big)=1 $ since the coupling preserves the order. %
So, thanks to Lemma~\ref{NeqEsp}, applied with \\
$\mathcal P = P_\L^- \opreceq P_\L^+ ( \ . \ | %
(\o^1(0),\o^2(0)) =(\minimalConf,\maximalConf) )$
and $Z= \o^2_0(n)-\o^1_0(n)$ we have: \\

$ P_\L^- \opreceq P_\L^+ \Big(  \o^1_0(n) \neq \o^2_0(n) \Big| %
(\o^1(0),\o^2(0)) =(\minimalConf,\maximalConf) \Big) %
 \leq  %
\constante \Big( %
P_\L^{\maximalConf}(\o_0(n) | \o(0)=\maximalConf ) - %
 P_\L^{\minimalConf}(\o_0(n) | \o(0)=\minimalConf ) \Big) $

where $\constante= (\min \{ s-s' : s > s'; s, s' \in \spin \})^{-1}$.
By Lemma~\ref{pptyCoup}, $\r(n)$ is bounded from
above by the l.h.s of the previous inequality.
We conclude by taking the limit in $n$, and using the finite volume
ergodicity.
\end{proof}
\begin{rem}
As an immediate consequence of Lemma~\ref{HWMErgodicite}
 we get $\lim_{n \to \infty} \r(n)=0$,
which implies the ergodicity of~$P$ thanks to Lemma~\ref{RhoControlErgod}.
\end{rem}

Let us denote by
\mbox{$ %
 R=\max_{k'\in V_0} \norme{k'} $}
the finite range of the local translation invariant PCA dynamics~$P$.
\begin{lemma} \label{InegRecurrence}  %
The following two inequalities hold:
\be{InegRho}
 \forall n \in \N^*,\ \r(2n) \leq  (2nR+1)^d \r^2(n) \ ;
\ee
\be{InegRecGeneNouvelle}
\forall n,\forall L\in \mathbb N^*, \ \r (2n) \leq 2(2L+1)^d \r^2(n) + %
{2 \constante} \Big( \int \s_0 \ d\nu_\BL^{\maximalConf} -  %
 \int \s_0 \ d\nu_\BL^{\minimalConf} \Big) \ .
\ee
\end{lemma}
\begin{proof}
Let $n$ be a fixed integer.

\noindent {\bf Proof of inequality~(\ref{InegRho})} \\
Let
$  \nu_n^{\minimalConf,\maximalConf}(\ .\ )=%
\coupling{(\o^1,\o^2)(n)\in .\ }{(\o^1,\o^2)(0)=%
(\minimalConf,\maximalConf)} $.
Using Markov property of $\coup{}$:
$$ %
\r(2n) = %
  \int \coupling{\o_0^1(2n) \neq \o_0^2(2n)}%
{(\o^1,\o^2)(n)=(\xi^{\minimalConf},\xi^{\maximalConf})} %
\nu_n^{\minimalConf,\maximalConf}(d\xi^{\minimalConf},d\xi^{\maximalConf}) %
\ . $$
Note that $\nu_n^{\minimalConf,\maximalConf}$-almost surely,
$\xi^{\minimalConf} \preceq %
\xi^{\maximalConf}$. %
Let $A = \{ (\xi^{\minimalConf},\xi^{\maximalConf} ) %
\tq \exists k \in \Zd ,\  \norme{k} \leq nR , %
 \ \xi^{\minimalConf}_k \neq \xi^{\maximalConf}_k \} $.
Thanks to Remark~\ref{discreteTime}
observe that the exact control of interaction information's
propagation for PCA implies that the above integral vanishes on $A^c$  because
\mbox{$\BnR \supset \overline{\{0\}}^{(n)}$}, and %
so \mbox{$ \xi^{\minimalConf}_{\BnR} \equiv \xi^{\maximalConf}_{\BnR}$}.
Then:
$$\r (2n)= \int_A \coupling{\o_0^1(n) \neq \o_0^2(n)}%
{(\o^1,\o^2)(0)=(\xi^{\minimalConf},\xi^{\maximalConf})}%
 \ \nu_n^{\minimalConf,\maximalConf}%
(d\xi^{\minimalConf},d\xi^{\maximalConf}) \ .$$
Using Lemma~\ref{pptyCoup}, we obtain
\mbox{$ \r (2n) \leq \r (n) \ \nu_n^{\minimalConf,\maximalConf}(A)$}. \\
Writing
\mbox{$ A= \cup_{\{k \in \Zd \tq  \norme{k} \leq nR \} } %
\{ (\xi^{\minimalConf},\xi^{\maximalConf}) \tq %
\xi^{\minimalConf}_k \neq \xi^{\maximalConf}_k \}$} %
we deduce:
$$\nu_n^{\minimalConf,\maximalConf} (A) \leq %
\sum_{k \in \Zd, \norme{k} \leq nR} %
\coupling{\o^1_k(n) \neq \o^2_k(n)}{(\o^1,\o^2)(0)=%
(\minimalConf,\maximalConf)}. $$
Since $P$ is translation invariant, the
conclusion follows from
$ \nu_n^{\minimalConf,\maximalConf} (A) 
  \leq   \r (n) \card{\BnR} \leq \r(n) (2nR+1)^d $
where $\card{\BnR}$ denotes the
cardinality of $\BnR$.
\medskip

\noindent {\bf Proof of inequality~(\ref{InegRecGeneNouvelle})} \\
Write
\mbox{$\r (2n)= \int \coupling{\o^1_0(2n) \neq \o^3_0(2n)}{%
(\o^1,\o^2,\o^3)(0)=(\minimalConf,\e,\maximalConf)}  \nu(d\e)$}
where $\nu$ is a $P$-stationary measure.
Note that %
\mbox{$ \o^1_0(n) \leq  \o^2_0(n) \leq \o^3_0(n) $}, \\
\mbox{$\coup{} \Big( %
(\o^1,\o^2,\o^3)\in \ . \ \Big| \ (\o^1,\o^2,\o^3)(0)= %
(\minimalConf,\e,\maximalConf) \Big)$}-almost surely,
so  that \\ %
\mbox{$ %
\{\o^1_0(n) \neq \o^3_0(n) \} =  %
\{ \o^1_0(n) \neq \o^2_0(n) \} \ \cup \  %
\{ \o^2_0(n) \neq \o^3_0(n) \} ,$}
 where the union is %
non necessarily disjoint
(unless cardinality of $\spin$ is $2$).
Thus, following decomposition holds:
\begin{eqnarray} \label{EqDecompPartielle}
{ \r (2n) \leq \int \coupling{\o^1_0(2n) \neq \o^2_0(2n)}%
{(\o^1,\o^2)(0)=(\minimalConf,\e)} \ \nu(d\e)} \nonumber  \\
\hspace{3cm}  +   \int \coupling{\o^1_0(2n) \neq \o^2_0(2n)}%
{(\o^1,\o^2)(0)=(\e,\maximalConf) }%
\ \nu(d\e)  \ . %
\end{eqnarray}

It is then enough to prove that each of these quantities are
bounded from above by half the quantity wanted.
Consider first the second term in the r.h.s. . \\
Let $
\nu_n^{\e,\maximalConf} = \coup{} \Big( %
(\o^1,\o^2)(n)=\ . \ \Big| \ %
(\o^1,\o^2)(0)=(\e,\maximalConf) \Big)$.
Let us write:
\begin{eqnarray*}
\lefteqn{ \int \coupling{\o^1_0(2n) \neq \o^2_0(2n)}%
{(\o^1,\o^2)(0)=(\e,\maximalConf) }%
 \ \nu(d\e) } 
\\
& =&   \iint \coupling{\o^1_0(n) \neq \o^2_0(n)}{ %
(\o^1,\o^2)(0)=(\xi^1,\xi^2)} %
\ \nu_n^{\e,\maximalConf}(d\xi^1,d\xi^2) \ \nu(d\e) \ .
\end{eqnarray*}
Let $L \in \N^*$ and \mbox{$A_L=\{(\xi^1,\xi^2) \in (\cvinf)^2 \tq %
({\xi^1})_\BL \equiv ({\xi^2})_\BL \} $}.
Let decompose the integration with respect to
$(\xi^1,\xi^2)$ into an
integration on $A_{L}^c$ and $A_{L}$. %
We will prove that:
\begin{eqnarray}
 (I) & = & \iint_{A_L^c} \coupling{\o^1_0(n) \neq \o^2_0(n)}{ %
(\o^1,\o^2)(0)=(\xi^1,\xi^2)} %
\ \nu_n^{\e,\maximalConf}(d\xi^1,d\xi^2)\ %
\nu(d\e) %
\nonumber \\
 & & \hspace{5cm} \leq  (2L+1)^d \r^2(n), 
\label{PartI}
\end{eqnarray}
\begin{eqnarray}
 (I\!I) & = &  \iint_{A_L} \coupling{\o^1_0(n) \neq \o^2_0(n)}{ %
(\o^1,\o^2)(0)=(\xi^1,\xi^2)} %
\ \nu_n^{\e,\maximalConf}(d\xi^1,d\xi^2) \ \nu(d\e) \nonumber \\
& & \hspace{5cm} \leq %
 \constante %
\Big( \int \s_0 \ d\nu_\BL^{\maximalConf} %
-   \int \s_0 \ d\nu_\BL^{\minimalConf} \Big)\hspace{-1mm}. %
 \label{PartII}
\end{eqnarray}
Let us consider part~$(I)$. Thanks to %
\mbox{$\nu_n^{\e,\maximalConf} (\xi^1 \preceq \xi^2)=1 $}
and using Lemma~\ref{pptyCoup}, we have \\
\mbox{$ (I)  \leq \r (n) \int \nu_n^{\e,\maximalConf} (A_L^c) \ \nu(d\e)$}.
Note that $ A_L^c$ may also be written
\mbox{$\cup_{k \in \BL} %
\{ (\xi^1,\xi^2) \tq ({\xi^1})_k \neq  ({\xi^1})_k \} $}.
Thus we have:
$$\nu_n^{\e,\maximalConf} (A_L^c) \leq \sum_{k \in \BL} %
\nu_n^{\e,\maximalConf} %
\{ (\xi^1,\xi^2) \tq ({\xi^1})_k \neq  ({\xi^2})_k \}\ .$$
Using translation invariance of the coupling and
 Lemma~\ref{pptyCoup}, the previous general term is equal to
$\coupling{\o^1_k(n) \neq \o^2_k(n)}{(\o^1,\o^2)(0)=(\e,\maximalConf)}  %
  \leq   \r (n) .$
So  \mbox{$\nu_n^{\e,+} (A_L^c) \leq \card{\BL} \ \r (n) $}, and then
(\ref{PartI}) follows.
\medskip

Part~$(I\!I)$: let $\t \in \cvinfLx{\BL}$ be fixed, and define
$ A_{L,\t} =\{ (\xi^1,\xi^2) \tq ({\xi^1})_\BL \equiv ({\xi^2})_\BL %
\equiv \t \} $.
So \mbox{$A_L=\miniop{}{\bigsqcup }{\t \in \cvinfLx{\BL}} A_{L,\t}$}
and following decomposition holds:
\begin{equation} (I\!I) =\int \sum_{\t \in \cvinfLx{\BL}} %
\int %
\coup{} \Big( \o^1_0(n) \neq \o^2_0(n) \Big| %
 (\o^1,\o^2)(0)=(\xi_1,\xi_2) \Big) %
 \ind{A_{L,\t}}{(\xi^1,\xi^2)}
 \nu_n^{\e,\maximalConf} (d\xi^1,d\xi^2) \ \nu(d\e) . \label{Principale}
\end{equation}
Let us now use the finite volume dynamics.
$\nu_n^{\e,\maximalConf}$ almost surely, we have %
\mbox{$\xi^1 \preceq \xi^2 $},
$({\xi^1})_\BL=({\xi^2})_\BL=\t$ %
and also
$\xi^2=\t ({\xi^2})_{\BL^c} \preceq \t (\maximalConf)_{\BL^c} $,
$ \t (\minimalConf)_{\BL^c} \preceq \xi^1=%
\t ({\xi^1})_{\BL^c} $.
Then:
\begin{eqnarray*}  P_\BL^{\minimalConf} %
\opreceq P \opreceq P \opreceq P_\BL^{\maximalConf} %
 ( %
\ \o^1 \preceq \o^2 \preceq \o^3 \preceq \o^4 \ \Big|  
(\o_\BL^1,\o^2,\o^3,\o_\BL^4)(0)= \\ \hspace{3cm} (\t, \t ({\xi^1})_{\BL^c}, \t%
({\xi^2})_{\BL^c}, \t ))=1
\end{eqnarray*}
which implies: 
\begin{eqnarray}
\lefteqn{ \coupling{\o^1_0(n) \neq \o^2_0(n)}{(\o^1,\o^2)(0)=%
(\t {\xi^1_{\BL^c}},\t {\xi^2_{\BL^c}})}} \nonumber \\
&  \leq & P_\BL^{\minimalConf} \opreceq P_\BL^{\maximalConf} %
 \big( {\o^1_0(n) \neq \o^2_0(n)%
 } \ | \ {(\o^1,\o^2)(0)=(\t,\t)} \big) \nonumber \\
& \leq & \constante  %
 \Big( %
P_\BL^{\maximalConf} ( \o_0(n) \ | \ \o_\BL(0)=\t ) - %
P_\BL^{\minimalConf} ( \o_0(n) \ | \ \o_\BL(0)=\t ) \Big)  , \label{Report1}
\end{eqnarray}
where the last inequality comes from
Lemma~\ref{NeqEsp} and from
the fact that \\
\mbox{$P_\BL^{\minimalConf} \opreceq P_\BL^{\maximalConf} %
 \Big(  \ . \  \Big| \ %
(\o^1,\o^2)(0)=(\t,\t) \Big)$}-almost surely, we have
\mbox{$\o^1_0(n) \leq  \o^2_0(n)$}.

On the other hand, note the following inequality:
\begin{eqnarray} \nu_n^{\e,\maximalConf} (A_{L,\t}) %
& = & \coupling{\o^1(n)_\BL \equiv \o^2_\BL(n) \equiv \t}%
{(\o^1,\o^2)(0)=(\e,\maximalConf)} \nonumber \\
&  \leq & \nu_n^{\e,\maximalConf} \Big(  (\xi^1,\xi^2) \tq %
({\xi^1})_\BL \equiv \t \Big) %
= P(\o_\BL(n)=\t | \ \o_\BL(0)=\e) \ . \label{Report2}
\end{eqnarray}
Reporting (\ref{Report1}) and (\ref{Report2}) in (\ref{Principale})
we find
\begin{eqnarray*}
 (I\!I) \leq \constante \int \sum_{\t \in \cvinfLx{\BL}} & &   %
 \hspace{-0,5cm} \Big( %
P_\BL^{\maximalConf} %
( \o_0(n) \ | \ \o_\BL(0)=\t ) - P_\BL^{\minimalConf} %
 ( \o_0(n) \ | \ \o_\BL(0)=\t ) %
\Big)  \\
& & \hspace{2cm} \ P(\o_\BL(n)=\t | \ \o_\BL(0)=\e)  \ \nu(d\e) \leq \constante \big( (a) - (b) \big) \ .
\end{eqnarray*}
We remark that %
\mbox{$ (a) =\int P \Big( f_{n,\maximalConf}(\o_\BL(n))\ %
\Big| \ \o_\BL(0)=\e \Big) %
\ \nu(d\e)  $}
with  \\ \mbox{$f_{n,\maximalConf}(\t) = %
P_\BL^{\maximalConf} ( \o_0(n) \ |\ \o_\BL(0)=\t )$}.
Using the fact that the function \mbox{$f_{n,\maximalConf}(.)$}
is increasing,  and Lemma~\ref{pptyCoup}
we state:
$$ (a) \leq \int \sum_{\t \in \cvinfLx{\BL}} %
P_\BL^{\maximalConf} ( \o_0(n) \ | \ \o_\BL(0)=\t ) \ %
 P_\BL^{\maximalConf} %
(\o_\BL(n)=\t \ |\ \o_\BL(0)=\e_\BL) \ \nu(d\e)  \ . $$
Using Markov property for the \mbox{$P_\BL^{\maximalConf}$} %
finite volume dynamics, we find:
\mbox{$ (a) \leq %
{\nu}(f_{2n,\maximalConf}) $}.
The function $f_{2n,\maximalConf}$ is increasing; thanks to
inequality~(\ref{IntermedEncadrementNu}), we
thus have
\mbox{$ (a) \leq \nu_\BL^{\maximalConf}(f_{2n,\maximalConf}) $}.
We can now write:
$$ (a) \leq \int P_\BL^{\maximalConf} %
 ( \o_0(2n) | \o_\BL(0)=\e_\BL ) \ \nu_\BL^{\maximalConf} (d\e_\BL) %
 =\int \s_0 \ %
d \nu_\BL^{\maximalConf} ,$$
where the last equality comes from the stationarity of $\nu_\BL^{\maximalConf}$
with respect to ${P_\BL^{\maximalConf}}$.
\medskip

Analogously we prove %
\mbox{$ (b) \geq %
\int \s_0 \ %
d \nu_\BL^{\minimalConf} $}.
Thus, the following inequality holds:
$$ (I\!I) \leq \constante \big( (a) - (b) \big) \leq %
\constante \left( \int \s_0 \ d \nu_\BL^{\maximalConf} %
- \int \s_0 \ %
d \nu_\BL^{\minimalConf} \right) ,$$
which gives the estimate of the second term in %
inequality~(\ref{EqDecompPartielle}). The first term is treated in %
the same way.
So the recursive inequality~(\ref{InegRecGeneNouvelle}) is established.
\end{proof}

We now  state some general analytic lemmas; for proofs see
~\cite{ThesePyl,MartinelliOlivieriI}.
\begin{lemma} \label{DecrescenzaPolinomiale} %
If $\lim_{n \to \infty} \r (n)=0$ and if
\mbox{$\exists \ (\tilde{C},M) \in (\mathbb R^+_*)^2,\ %
\exists L_1 \in \mathbb N^*, \forall L \in \mathbb N^*, L \geq L_1,  \forall n \in \N^*$}
$$\r(2n) \leq 2(2L+1)^d \r (n)^2 + 2\tilde{C}e^{-ML}  $$
then \mbox{$\lim_{n \to \infty} n^d \r(n)=0$}.
\end{lemma}
\begin{lemma} \label{VitesseExplicite}
If
\mbox{$\lim_{n \to \infty} n^d \r(n)=0$}, and if
inequality~(\ref{InegRho}) holds
then, for all $n_1$ such that
\mbox{$(2^d\hat{C})\ n_1^d \r(n_1)<1$},
we have:
$$ \forall n \geq n_1,\ \r(n) \leq e^{-\lambda n}$$
where \mbox{$\lambda=-\frac{1}{2n_1} \log (2^d \hat{C} n_1^d \r(n_1)) >0$}.
\end{lemma}

\section{\label{SectionApplication}\label{SectionProof2}Proof of the Theorem~\ref{ErgodClasseReversible}}

For general PCA
in finite volume, invariant
measures are not explicitly known; but for the class \nosPCA here considered,
we computed them ({\it cf.}  Proposition~3.1 {\it in}~\cite{DPLR}).
The unique reversible measure for the PCA dynamics $P_\L^\t$ is defined by
\be{DefNuVolFin}
 \nu_{\L}^{\t}(\s) = \frac{1}{\W_{\L}^{\tau}} \prod_{k \in \L} %
 \cosh \left( \b \sum_{j \in \Z^{d}} \K(k-j)
 \tilde{\s}_{j} \right) e^{\b \s_{k} \sum_{j \in \L^{c}} \K(k-j)
 \tau_{j}},
\ee
 where \mbox{$\tilde{\s}=\s_{\L} \t_{\L^c}$},
 and $\W_{\L}^{\t}$ is the normalisation factor.
Such measure
 does not coincide with the finite volume Gibbs measures
$\mu_\L^\t(\s)=\frac{1}{\Zp_\L^\t} \textrm{exp} ( - \sum_{A \subset %
{\mathbb Z^d}, A \cap \L \neq \emptyset} \varphi_A(\s_\L \t_{\L^c}))$
contrary to what happens for Glauber dynamics when detailed balance
holds.
Nevertheless, they are related as relation~(\ref{GibbsStatVolFin}) attempts.
We will not write down all technical computations
which prove relations
  (\ref{GibbsStatVolFin}), %
 (\ref{CondPlus}). %
Interested reader may refer respectively  to
Proposition~4.1.8 %
and Property~4.1.12 {\it in}~\cite{ThesePyl}.
\medskip

Let
 $\L,\L'$ two finite subsets of $\Zd$ such that $\L \subset \L'$ and
$\partial_i \L \cap \partial_i \L' = \varnothing$, %
where
$\di{\L} \eqdef \{ k \in \L \tq U_k \cap {\L}^c \neq \emptyset \}$.
 Let $\t'$ be a boundary condition of $\L$  and  $\mu_{\L}^{\t'}$ denotes the finite volume Gibbs
distribution associated to the potential~$\varphi$
on the volume $\L$ with boundary condition $\t'$. We then state:
\be{GibbsStatVolFin}
\nu_{\L'}^\t (d\s_\L \vert \s_{\L' \setminus \L} ) = %
\mu_\L^{\s_{\L' \setminus \L}\t_{\L'^c}} (d \s_\L) \ .
\ee
Note that the potential $\varphi$ is not
really a ferromagnetic potential in the usual sense. However we can
check that associated finite
volume Gibbs
measures verify a kind of monotone behaviour:
$ \t_1 \preceq \t_2 \Rightarrow %
\mu_\L^{\t_1} \preceq \mu_\L^{\t_2}$ (see Proposition~4.1.9 {\it in}~\cite{ThesePyl}).
In particular, Gibbs measures on~$\cvinf$ obtained
as $\mu^{+}= \lim_{\L \nearrow \Zd} \mu_\L^{(+)_{\L^c}}$
and $\mu^{-}=\lim_{\L \nearrow \Zd} \mu_\L^{(-)_{\L^c}}$
are extremal states in the sense of stochastic ordering of the set $\G{\varphi}$.
Recall $\mu$ probability measure on $\cvinf$ is in $\G{\varphi}$ if, %
{\it per definitionem}, for any finite volume $\L \subset \Zd$,
 a version of the conditioned measure $\mu (d\s_\L |\s_{\L^c})$ is $\mu_\L^{\s_{\L^c}}(d\s_\L)$.
Finally, let us state the following lemma:
\begin{lemma} \label{ImplicationLimites}
If the Weak Mixing Condition~\WM holds for the %
potential~$\varphi$ associated to the PCA dynamics~$P$,
then assumption \assump holds %
for $P$.
\end{lemma}
\begin{proof}
According to the finite range $R$, %
let $L>R$.
It is enough to show \\
$ 
\Big( \int \s_0 \ d\nu_\BL^+ -   \int \s_0 \ d\nu_\BL^- \Big) \leq %
\Big( \int \s_0 \ d\mu_\BLR^+ -   \int \s_0 \ d\mu_\BLR^- \Big) $.
Let us first check \\ \mbox{$ %
 \int \s_0 \ d\nu_\BL^+ \leq \int \s_0 \ d\mu_\BLR^+$}.
Let $f_0$ be the increasing function defined on~$\cvinf$ by %
$f_0(\s)=\s_0$.
Note
\mbox{$\int \s_0 \ d \nu_\BL^+ = %
\nu_\BL^+( \ \nu_\BL^+ (\ f_0 \ | \s_{\BL \setminus \BLR} )) $}.
Using relation~(\ref{GibbsStatVolFin}) with $\L'=\BL$ and $\L=\BLR$, we then get
$ \nu_\BL^+(f_0)=\nu_\BL^+( %
\mu_{\BLR}^{ \s_{\BL \setminus \BLR}(\+)_{\BL^c}}%
(f_0)) $.
On the other hand, using the monotonicity in the
 boundary condition of the finite
volume Gibbs measures, we find
$  \mu_{\BLR}^{ \s_{\BL \setminus \BLR}(\+)_{\BL^c}}%
(f_0) \leq \mu_{\BLR}^{(\+)_{\BLR^c}}(f_0)$.
So desired inequality holds.
\mbox{$ \nu_\BL^-(f_0) \geq  \mu_\BLR^-(f_0) $} can be analogously checked.
\end{proof}
\begin{lemma} \label{CoincidenceMesuresExtremes}
For a PCA dynamics $P$ of class \nosPCA with $\K(.)$ non negative,
the extremal stationary measures $\nuinf,\nusup$ coincide respectively with
extremal Gibbs measures $\mu^-$ and $\mu^+$ of $\G{\varphi}$
(possibly these four measures coincide).
\end{lemma}
\begin{proof}
Let $\L,\ \L'$ be two finite subsets of $\Zd$ such that
$\L \subset \L'$. Then, for all configurations
\mbox{$\s_{\L' \setminus \L} \in \cvinfLx{\L'\setminus \L}$}, %
finite volume reversible measures with extremal boundary condition
are such that:
\begin{equation} \label{CondPlus}
 \nu_{\L'}^+ \left( (.)_\L \vert \s_{\L' \setminus \L} \right) %
\preceq \nu_\L^+(.) \ ; \quad %
 \nu_{\L'}^- \left( (.)_\L \vert \s_{\L' \setminus \L} \right) %
\succeq \nu_\L^-(.)
\end{equation}
(see Property~4.1.12 {\it in}~\cite{ThesePyl} for a precise proof).
Using relation~(\ref{GibbsStatVolFin}), we can deduce
from the previous result the following inequalities
between finite volume Gibbs measure and reversible measure,
with extremal boundary condition:
\mbox{$ \mu_\L^+ \preceq \nu_\L^+$} and
\mbox{$ \mu_\L^- \succeq \nu_\L^-$}.
Taking now the limit in volume, we find:
\mbox{$ \mu^+ \preceq \nu^+$} and
\mbox{$ \mu^- \succeq \nu^-$}.
\medskip

On the other hand, $\nu_\L^+$ is $P_\L^+$-reversible, so taking the limit,
$\nu^+$ is $P$-reversible. Analogously, $\nu^-$ is $P$-reversible.
From~$\R=\S \cap \G{\varphi}$,  we conclude $\nu^-$ and $\nu^+$ are
Gibbs measures, so thanks to the fact that $\mu^-$ and $\mu^+$ are %
stochastic ordering extremal states for Gibbs measures, we deduce:
  \mbox{$ \nu^+ \preceq \mu^+$} and
\mbox{$ \mu^- \preceq \nu^-$}. Thus the conclusion follows.
\end{proof}

Here is the proof of Theorem~\ref{ErgodClasseReversible}:

\begin{proof}
When there is phase transition, since
$\mu^-$ and $\mu^+$ are extremal states
for $\G{\varphi}$, we have that $\mu^- \neq \mu^+$. So, using
Lemma~\ref{CoincidenceMesuresExtremes}, the two reversible (also stationary)
measures $\nu^-$ and $\nu^+$ are different. Then, dynamics $P$ can not be ergodic.
\medskip

When there is no phase transition, then \mbox{$\G{\varphi}=\{\mu\}$}
where \mbox{$\mu=\mu^-=\mu^+$} is the unique Gibbs state.
Thanks to~Lemma~\ref{CoincidenceMesuresExtremes}, it holds $\nu^-=\mu^-=\mu^+=\nu^+$.
The Proposition~\ref{ResultatsDeMonotonie}
states the uniqueness of the $P$-stationary measure and
the ergodicity of the PCA dynamics~$P$.
\medskip

Finally, if weak mixing condition~\WM
is assumed, then Lemma~\ref{ImplicationLimites} implies that inequality~\assump holds. We conclude
using Theorem~\ref{ergodicite}.
\end{proof}


{\bf \large Acknowledgments} :

This work is part of the author's PhD Thesis, realized at the
Universit\'e Lille~1 and Politecnico of Milan.
P.-Y. Louis warmly thanks
his PhD advisers, P.~Dai~Pra and S.~R\oe lly, for supervising his
work, and for the encouragements they provided. An anonymous referee
is acknowledged for carefully reading the first version of this paper.

The Mathematics' Department of Padova University and the Interacting Random Systems
group of Weierstrass Institute for Applied Analysis and
Stochastics in Berlin, where part
of this work was done, are acknowledged for their kind hospitality.


\bibliographystyle{amsalpha}


\vspace*{5mm}

\end{document}